\documentclass{article}

\usepackage{amsmath, amsthm, amssymb}

\newcommand{\mc}{\mathcal}
\newcommand{\sps}{\left[\begin{array}{l}p\\s\end{array}\right]}
\newcommand{\hoz}{(H^{'}_{p-1})_0}
\newcommand{\hou}{(H^{'}_{p-1})_1}
\newcommand{\mpd}{\;\;\text{mod}\,p^2}
\newcommand{\mpu}{\;\;\text{mod}\,p}

\author{Claire Levaillant}
\title{Powers of two weighted sum of the first $p$ divided Bernoulli numbers modulo $p$}

\begin{document}
\maketitle
\begin{center}
\underline{Abstract}
\end{center}
We show that, modulo some odd prime $p$, the powers of two weighted sum of the first $p-2$ divided Bernoulli numbers equals the Agoh-Giuga quotient plus twice the number of permutations on $p-2$ letters with an even number of ascents and distinct from the identity. We provide a combinatorial characterization of Wieferich primes, as well as of primes $p$ for which $p^2$ divides the Fermat quotient $q_p(2)$.

\section{Introduction}

Congruences involving Bernoulli numbers or divided Bernoulli numbers have drawn the attention of many mathematicians. The list is so long that we may not venture citing any author here. In this note, we study powers of two weighted sums of the first $p$ divided Bernoulli numbers modulo $p$, where $p$ is a given odd prime number. Bernoulli numbers are rational numbers and presently we view them as $p$-adic numbers. We define the divided Bernoulli numbers as
$$\mathcal{B}_n=\frac{B_n}{n}\in\mathbb{Q}_p\;\text{when $n\geq 1$ and}\;\;\mathcal{B}_0=-\frac{1}{p}\in\mathbb{Q}_p,$$
where $\mathbb{Q}_p$ denotes the field of $p$-adic numbers. The convention for $\mathcal{B}_0$ is our own convention, especially the minus sign which is simply there because it will allow for a nicer formula.
We recall that all the Bernoulli numbers with odd indices are zero, except $B_1$ which takes the value $-\frac{1}{2}$.

After excluding both the upper bound and the lower bound of the weighted sum, we are able to express the remaining sum in terms of the Agoh-Giuga quotient and of some Eulerian sum modulo $p$. The Agoh-Giuga quotient is the $p$-adic integer
$$\frac{p\,B_{p-1}+1}{p}\in\mathbb{Z}_p$$
Von Staudt-Clausen's theorem, a theorem dating from $1840$ that was independently proven by von Staudt and Clausen, asserts that the Bernoulli numbers $B_{2k}$ sum to zero when added all the fractions $\frac{1}{q}$ with $q$ prime such that $q-1$ divides $2k$.
It is a consequence of von Staudt-Clausen's theorem that $p\,B_{p-1}$ is a $p$-adic integer of residue $-1$ modulo $p\,\mathbb{Z}_p$. \\
Eulerian sums are sums involving Eulerian numbers. The Eulerian number $E(m,n)$ is by definition the number of permutations on $n$ letters with $m$ ascents. If $(i_1,i_2,\dots,i_n)$ denotes the permutation mapping the integer $j$ onto $i_j$, an ascent is when $i_{k+1}>i_k$.
Eulerian numbers originate in Euler's book dating from $1755$ in which he investigates shifted forms of what are now called Eulerian polynomials:
$$E_n(t)=\sum_{m=0}^n E(n,m)\,t^m$$
They get more attention a century later in the studies of Worpitzky, see \cite{WO2}.
\\

Our main result is the following.
\newtheorem*{Main Theorem}{Main Theorem}
\begin{Main Theorem}
Let the $\mathcal{B}_k$'s denote the divided Bernoulli numbers, as defined earlier. Then, we have
\begin{equation}\sum_{k=1}^{p-2}\frac{\mc{B}_k}{2^k}=\frac{p\,B_{p-1}+1}{p}+2\left(\begin{array}{l}
\#\text{permutations on $(p-2)$ letters}\\\text{with an even number of ascents}\\\text{and distinct from the identity}
\end{array}\right)\;\text{mod}\,p\end{equation}
\end{Main Theorem}
\noindent The second component of the sum to the right hand side of $(1)$ is precisely twice the Eulerian sum
$$\sum_{k=0}^{\frac{p-5}{2}} E(p-2,2k)$$
Eulerian numbers can be computed from the recursion formula
$$E(n,m)=(n-m)\,E(n-1,m-1)+(m+1)\,E(n-1,m)$$
We prove more results including but not limited to the Fermat quotient $q_p(2)$ to the moduli $p$ and $p^2$ and congruences concerning sums of powers of odd integers and sums of divided Bernoulli numbers.

Some of the proofs in the current paper are based on using two equivalent definitions for the unsigned Stirling numbers of the first kind. From these two definitions we can derive the following lemma.
\newtheorem{Lemma}{Lemma}
\begin{Lemma}
\begin{eqnarray}
p!&=&\sum_{s=1}^p\left[\begin{array}{l} p\\s\end{array}\right]\\
(p+1)!&=&\sum_{s=1}^p\left[\begin{array}{l} p\\s\end{array}\right]2^s
\end{eqnarray}
\end{Lemma}
\noindent Equality $(2)$ comes from the fact that $\sps$ counts the number of permutations of $p$ objects that decompose into a product of $s$ disjoint cycles. Since any permutation of $p$ objects can be uniquely written as a product of cycles with disjoint supports, equality $(2)$ holds. \\
As for Equality $(3)$, it arises from the other definition of the unsigned Stirling numbers of the first kind, namely, $\sps$ is the unsigned coefficient of $x^s$ in the falling factorial
$$x(x-1)\dots(x-(p-1))$$ and so, one can write:
\begin{equation}
x(x-1)(x-2)\dots (x-(p-1))=\sum_{s=1}^p (-1)^{s-1}\sps\,x^s\end{equation}
It now suffices to specialize $x=-2$ in Equality $(4)$. We thus get Equality $(3)$ of Lemma $1$.
Modulo $p^2$, the two factorials $p!$ and $(p+1)!$ of the lemma are equal.

By \cite{LEV2} or \cite{GL}, we know how to calculate the Stirling numbers modulo $p^2$ in terms of Bernoulli numbers. Then by using some results of Emma Lehmer \cite{LE}, we are able to relate the investigated weighted sum to a sum of powers of integers modulo $p^2$, just as stated in the following proposition.
\newtheorem{Proposition}{Proposition}
\begin{Proposition}
\begin{equation}\sum_{a=1}^{\frac{p-1}{2}}a^{p-1}=-1-\frac{1}{2}\,p\,B_{p-1}+\frac{p}{2}\Bigg(\frac{3}{2}+\sum_{k=1}^{\frac{p-3}{2}}\frac{B_{2k}}{k\,2^{2k}}\Bigg)\;\text{mod}\,p^2
\end{equation}
\end{Proposition}

Independently, by working out the sum to the left hand side of $(5)$ modulo $p^2$ and by using some known formula relating an alternating sum of Eulerian numbers to the Bernoulli numbers, we obtain an intermediate result on the Fermat quotient
$$q_p(2)=\frac{2^{p-1}-1}{p}$$
The result is stated in the following theorem.
\newtheorem{Theorem}{Theorem}
\begin{Theorem}
\begin{equation}
q_p(2)=2\,\left(\begin{array}{l}\#\text{permutations on $(p-2)$ letters}\\\text{with an even number of ascents}\\\text{and distinct from the identity}\end{array}\right)+1\;\;\text{mod}\,p
\end{equation}
\end{Theorem}
Forthcoming Theorem $9$ of the introduction provides the next term in the $p$-adic expansion of the Fermat quotient $q_p(2)$, expressed combinatorially and in terms of the residue of the Agoh-Giuga quotient.

As a corollary obtained from Eisenstein's formula \cite{EI} for the Fermat quotient
\begin{equation*}
q_p(2)=\frac{1}{2}\sum_{k=1}^{p-1}\frac{(-1)^{k-1}}{k}\;\text{mod}\,p,
\end{equation*}

we have:
\begin{Theorem}
\begin{equation}
\sum_{k=1}^{p-1}\frac{(-1)^k}{k}=2\Bigg(1-2\left(\begin{array}{l}\#\text{permutations on $(p-2)$ letters}\\\text{with an even number of ascents}\end{array}\right)\Bigg)\;\;\text{mod}\,p
\end{equation}
\end{Theorem}
\noindent And since we have $H_{p-1}=0\;\text{mod}\,p$ by Wolstenholme's theorem \cite{WO1}, where $H$ denotes the harmonic number, we thus deduce
\begin{Theorem}
$$q_p(2)=\sum_{\begin{array}{l}k=1\\\text{k odd}\end{array}}^{p-1}\frac{1}{k}=2\,\left(\begin{array}{l}\#\text{permutations on $(p-2)$ letters}\\\text{with an even number of ascents}\\\text{and distinct from the identity}\end{array}\right)+1\;\;\text{mod}\,p$$
\end{Theorem}

\newtheorem{Remark}{Remark}
\begin{Remark}
The harmonic sum of Theorem $3$ with odd reciprocals is sometimes denoted by $H^{'}_{p-1}$ in the literature. Back in $1900$, Glaisher has obtained the congruence to the left in Theorem $3$ to the modulus $p^2$. Namely, he showed in \cite{GL2} that
$$q_p(2)=H^{'}_{p-1}+\frac{1}{2}\,p\,(H^{'}_{p-1})^2\;\;\text{mod}\,p^2$$
Vandiver provided a simpler proof for the latter congruence in \cite{VA}. \\At the end of his Remark $5.3$ of \cite{SU}, Sun provides a congruence for $H^{'}_{p-1}$ to the modulus $p^3$. It is the following:
$$H^{'}_{p-1}=q_p(2)-\frac{1}{2}\,pq_p^2(2)+\frac{1}{3}p^2q_p^3(2)-\frac{1}{24}\,p^2\,B_{p-3}\;\;\text{mod}\,p^3$$

\end{Remark}

\noindent Further, from Stern's formula \cite{ST} modulo $p$ for the Fermat quotient $q_p(2)$, namely
$$q_p(2)=\sum_{k=1}^{\frac{p-1}{2}}\frac{(-1)^{k-1}}{k}\;\text{mod}\,p,$$
\noindent joint with Lehmer's formula for the harmonic number $H_{\frac{p-1}{2}}$ modulo $p^2$,
we derive
\begin{Theorem}
$$\sum_{\begin{array}{l}k=1\\\text{k odd}\end{array}}^{\frac{p-1}{2}}\frac{1}{k}=\frac{1}{2}\Bigg(1-2\,\left(\begin{array}{l}\#\text{permutations on $(p-2)$ letters}\\\text{with an even number of ascents}\end{array}\right)\Bigg)\;\text{mod}\,p$$
\end{Theorem}

Coming back to the sketch of the proof of the main theorem, 
from Eq. $(3)$ of Lemma $1$ written modulo $p^2$ and Theorem $3$, we derive the formula of the main theorem.
Stated under the form of the title of the current paper, it reads instead:
\begin{Theorem} Under the convention that $\mathcal{B}_0=-\frac{1}{p}$, we have:
\begin{equation}
\sum_{k=0}^{p-1}\frac{\mathcal{B}_k}{2^k}=0\;\mpu
\end{equation}
\end{Theorem}
Last, a consequence of Eq. $(2)$ from Lemma $1$ is that the sum of the first $(p-2)$ divided Bernoulli numbers is congruent to the residue of the Agoh-Giuga quotient modulo $p$ minus $1$, that is congruent to the residue of the Wilson quotient. We also gather more congruences in the theorem below.
\begin{Theorem}Let $\mathcal{B}_i$ denote the $i$th divided Bernoulli number. Let $w_p$ denote the Wilson quotient and let $q_p(2)$ denote the Fermat quotient. Then, we have:
\begin{eqnarray}\sum_{k=1}^{p-2}\mathcal{B}_k&=&\frac{p\,B_{p-1}+1}{p}-1\;=\;w_p\qquad\qquad\;\;\,\,\,\text{mod}\,p\\
\sum_{k=1}^{p-2}\frac{\mathcal{B}_k}{2^k}&=& q_p(2)+w_p\;\;\qquad\qquad\qquad\qquad\;\;\;\;\,\,\text{mod}\,p\\
\sum_{k=1}^{p-2}\frac{\mathcal{B}_k}{2^k}&=&-\frac{1}{2}H_{\frac{p-1}{2}}+\frac{p\,B_{p-1}+1}{p}-1\;\;\;\,\qquad\text{mod}\,p\\
\sum_{\begin{array}{l}1\leq k\leq p-2\\\text{k is even}\end{array}}\frac{B_k}{k\,2^k}&=&\frac{3}{4}+\sum_{\begin{array}{l}1\leq k\leq p-2\\\text{k is even}\end{array}}\Big(\frac{B_k-1}{k}\Big)\;\;\;\;\text{mod}\,p
\end{eqnarray}
\end{Theorem}
\noindent In the statement above, $(10)$ and $(12)$ are obtained by combining several of the previous results. Whereas $(11)$ is a rewriting of $(10)$ using the fact that \begin{equation} H_{\frac{p-1}{2}}=-2q_p(2)\;\;\text{mod}\,p\end{equation}
This is for instance a special case of Vandiver's congruence dating from $1917$. His result of \cite{VA} is built out of Lerch's original formula dating from $1905$ and asserts that
$$n\,q_p(n)=-\sum_{j=1}^{n-1} H_{\lfloor \frac{jp}{n}\rfloor}\;\;\text{mod}\,p$$
We will see along the paper that Lehmer has also obtained Congruence $(13)$ to the next modulus $p^2$, and Sun in \cite{SU} even gets a congruence modulo $p^3$.
Regarding $(10)$, Beeger computed $w_p$ modulo $p$ for $p<300$. For this range of primes, $w_p=0\;\text{mod}\,p$ only when $p=5$ and $p=13$, see \cite{BE}. For these values of $p$, we see with $(10)$ that the powers of two weighted sum of divided Bernoulli numbers equals the Fermat quotient modulo $p$. In particular,
\begin{equation}
\sum_{k=1}^{11}\frac{B_k}{k2^k}=q_{13}(2)\;\;\text{mod}\,13
\end{equation}

Additionally, we prove the following statements concerning sums of powers.
\begin{Theorem}\hfill\\\\
$(i)$ The sum of the first $\frac{p-1}{2}$ odd powers of the first $\frac{p-1}{2}$ integers is congruent to $-\frac{1}{2}$ modulo $p$.\\
$(ii)$ The sum of the first $\frac{p-1}{2}$ odd powers of the next $\frac{p-1}{2}$ integers is congruent to $\frac{1}{2}$ modulo $p$.
\end{Theorem}
For the next statement, we must introduce a few notations.
\newtheorem{Notation}{Notation}
\begin{Notation}
Denote by $N_{p-2}$ the number of permutations on $p-2$ letters with an even number of ascents.
\end{Notation}
\begin{Notation}
Let $x$ be a $p$-adic integer. We denote by $(x)_k$ the $(k+1)$-th coefficient in the $p$-adic expansion of $x$, that is
$$x=\sum_{j=0}^{\infty}(x)_j\,p^j$$
In particular, $(x)_0$ is the residue of $x$ modulo $p$. We will make an exception with $(pB_{p-1})_1$ which will denote the Agoh-Giuga quotient.
\end{Notation}
The following theorem addresses sums of powers of odd integers modulo $p^2$.
\begin{Theorem}
\begin{eqnarray}(i)\;\; \sum_{\begin{array}{l}x=1\\\text{x odd}\end{array}}^{p-1} x^{p-2}&=& (2N_{p-2})_0-1+p\,\Bigg((p\,B_{p-1})_1+(2N_{p-2})_1-\big((2N_{p-2})_0-1\big)^2-2\Bigg)\notag\\&&\qquad\qquad\qquad\qquad\qquad\qquad\qquad\qquad\;\text{mod}\,p^2\\
(ii)\;\;\sum_{\begin{array}{l}x=1\\\text{x odd}\end{array}}^{p-1} x^{p-1}&=& -\frac{1}{2}+\frac{p}{2}\Bigg((p\,B_{p-1})_1-(2N_{p-2})_0+1\Bigg)\;\;\text{mod}\,p^2\end{eqnarray}
\end{Theorem}
By Beeger's result mentioned earlier, we know that $(13\,B_{12})_1=1$. \\From the table below listing the Eulerian numbers with an even number of ascents computed from the Euler triangle,
$$\begin{tabular}{|c|c|c|c|c|c|c|}
\hline
n\textbackslash m&0&2&4&6&8&10\\
\hline
11&1&152\,637&9\,738\,114&9\,738\,114&152\,637&1\\
\hline
\end{tabular}$$
we compute that $N_{11}=19\,781\,504$ and $(2N_{11})_0=4$ and $(2N_{11})_1=8$. \\
We can check that Congruence $(15)$ is verified when $p=13$. Namely, we have:
$$1^{11}+3^{11}+5^{11}+\dots+11^{11}=(2N_{11})_0-1+13\Big((2N_{11})_1-(2N_{11})_0^2+2(2N_{11})_0-2\Big)\;\;\text{mod}\,13^2$$

Finally, we confront our congruence $(15)$ of Theorem $8$ and a result of Emma Lehmer in order to obtain the next term in the $p$-adic expansion of the Fermat quotient $q_p(2)$. Our result is the following.

\begin{Theorem}
\begin{multline}
2^{p-1}=1+p\Bigg((2N_{p-2})_0-1\Bigg)\\+\;p^2\Bigg([(p\,B_{p-1})_1][(2N_{p-2})_0]-(2N_{p-2})_0^2+(2N_{p-2})_0+(2N_{p-2})_1-2\Bigg)\;\;\text{mod}\,p^3
\end{multline}
\end{Theorem}

\noindent Confronting our Theorem $9$ and Remark $1$, we obtain:
\begin{Theorem}
Let $H^{'}_{p-1}=1+\frac{1}{3}+\frac{1}{5}+\dots +\frac{1}{p-2}$. We have modulo $p$:
\begin{equation}
(2N_{p-2})_0=\hoz+1\end{equation}
\begin{equation}\begin{split}(2N_{p-2})_1+\big((2N_{p-2})_0-1)_1=&\frac{3}{2}\hoz^2+\hou+\hoz\\&-(p\,B_{p-1})_1\Big(1+\hoz\Big)+2\end{split}
\end{equation}

\end{Theorem}
By computing $(2N_{p-2})_0$ in another way, we obtain the following interesting congruence.
\begin{Theorem}
$$\sum_{\begin{array}{l}m=1\\\text{m odd}\end{array}}^{p-2}H_m=\frac{1}{2}\big(1+(H^{'}_{p-1})_0\big)=\frac{1+q_p(2)}{2}\mpu$$
\end{Theorem}

A prime $p$ is called a Wieferich prime if $q_p(2)=0\;\text{mod}\,p$. In other words, a Wieferich prime satisfies to Fermat's little theorem at the next modulus of $p$. Wieferich had introduced these primes in $1909$ while working on Fermat's last theorem. He had shown that if there exist solutions to $x^p+y^p=z^p$ in integers $x,y,z$ and $p$ is an odd prime with $p$ coprime to $x$, $y$, $z$, then $p$ is a Wieferich prime. As of September $2018$, the only known Wieferich primes are $1093$ (found by Meissner \cite{MI} in $1913$ and only such prime below $2000$) and $3511$ (found by Beeger \cite{BE} in $1922$). In \cite{DO}, the authors show that there are no other Wieferich primes $p<6.7\times 10^{15}$, yet improving on the bounds provided by a dozen of other mathematicians (the latest bound had been reached in $2005$). Previously, Grave had conjectured that Wieferich primes did not exist. In $1988$, Joseph Silverman showed in \cite{SI} that if the abc conjecture holds, then there exist infinitely many non-Wieferich primes.\\

We now present a combinatorial characterization of a Wieferich prime. In light of Theorem $3$, we have:
\begin{Theorem}
Let $N_{p-2}$ denote the number of permutations on $p-2$ letters with an even number of ascents. Let $(2N_{p-2})_0$ denote the residue modulo $p$ of twice this number. Then,
\begin{eqnarray}
\text{$p$ is a Wieferich prime}&\Leftrightarrow& (2N_{p-2})_0=1\\&\Leftrightarrow& \hoz=0\end{eqnarray}
\end{Theorem}
\newtheorem{Definition}{Definition}
\begin{Definition} (super Wieferich prime)
We will call a prime a "super Wieferich prime" if $q_p(2)=0\;\text{mod}\,p^2$
\end{Definition}
In light of Theorems $9$ and $10$, we have:
\begin{Theorem}
Let $N_{p-2}$ denote the number of permutations on $p-2$ letters with an even number of ascents. Let $(2N_{p-2})_0$ denote the residue modulo $p$ of twice this number. Let $(2N_{p-2})_1$ denote the next residue in the $p$-adic expansion of twice this number. Then,
\begin{eqnarray}\bigg(\text{$p$ is a super Wieferich prime}\bigg)&\Leftrightarrow&\left\lbrace\begin{array}{l}\!\!\!(2N_{p-2})_0=1\;\qquad\qquad\;\;\mpu\\
\!\!\!(2N_{p-2})_1+(p\,B_{p-1})_1=2\!\!\!\mpu\end{array}\right.\\
&&\notag\\
&\Leftrightarrow& \hoz=\hou=0
\end{eqnarray}
\end{Theorem}
We deduce a combinatorial necessary condition for a prime to be a super Wieferich prime and a Wilson prime.
\begin{Theorem} (Search for Wilson and super Wieferich primes).\\\\
Let $N_{p-2}$ denote the number of permutations on $p-2$ letters with an even number of ascents. Let $(2N_{p-2})_0$ denote the residue modulo $p$ of twice this number. Let $(2N_{p-2})_1$ denote the next residue in the $p$-adic expansion of twice this number.
\begin{center} If $p$ is a Wilson and super Wieferich prime, then $(2N_{p-2})_0=(2N_{p-2})_1=1$\end{center}
\end{Theorem}
We draw below a list of open questions.

\newtheorem{Open problem}{Open problem}
\begin{Open problem}
Does there exist any Wilson prime that is also a Wieferich prime ?
\end{Open problem}
\begin{Open problem}
Does there exist any super Wieferich prime ?
\end{Open problem}
If the answer to both questions is yes, we ask:
\begin{Open problem}
Does there exist any Wilson prime that is a super Wieferich prime ?
\end{Open problem}

In \cite{CDP}, Richard Crandall, Karl Dilcher and Carl Pomerance reported that there were no other Wilson primes than $5$, $13$ and $563$ less than $5\times 10^8$. \\The authors worked on both the Wieferich search and the Wilson search, but they did not study whether the two sets intersect or not. It has been conjectured that infinitely many Wilson primes exist, as for Wieferich primes. \\

The paper is structured as follows. \\

In $\S\,2$, we prove the main theorem as well as Theorems $1$-$6$ concerning the residues of the Fermat quotient $q_p(2)$ (Theorem $1$), resp the alternating harmonic sum and the harmonic sum of odd reciprocals corresponding to the harmonic number $H_{p-1}$ (Theorems $2$ and $3$), resp the harmonic sum of odd reciprocal corresponding to the harmonic number $H_{\frac{p-1}{2}}$ (Theorem $4$), resp the powers of two weighted sum of the first $p$ divided Bernoulli numbers (Theorem $5$), resp a sum of divided Bernoulli numbers whose indices range between $1$ and $p-2$ (Theorem $6$).

The next paragraph $\S\,3$ deals with sums of powers of integers exclusively. We prove Theorems $7$ and $8$.

In $\S\,4$, we prove Theorem $9$ which provides a congruence for the $p$-adic expansion of $q_p(2)$ up to the modulus $p^3$. We deduce Theorem $10$ by using Glaisher's expression for the Fermat quotient modulo $p^2$, taken from \cite{GL2}. The results of Theorems $12$ and $13$ can then be read out of Theorems $9$ and $10$. At the end of this part, we also prove Theorem $11$, providing a congruence relating modulo $p$ the Fermat quotient to a sum of harmonic numbers with odd indices. The result follows from Cong. $(18)$ of Theorem $10$ and from a separate calculation of the Eulerian numbers modulo $p$. The latter computation is the purpose of forthcoming Proposition $4$ in $\S\,4$.

In the last part, we mention Zhi-Wei Sun's work dealing with similar powers of two weighted sums, concerning divided harmonic numbers this time. We came across his work after this work was completed. \\ Tentatively connecting our work and Sun's work is left for future plans.\\
As part of his results, Z-W. Sun proves that the powers of two weighted sum of the first $(p-1)$ divided harmonic numbers is congruent to zero modulo $p$. Here we investigate when a similar statement holds if we consider a powers of two weighted sum of the first $(p-2)$ divided Bernoulli numbers. Further, we generalize this investigation to powers of integers $k$ with $3\leq k\leq p-1$. \\
Our conclusions are gathered in the following theorem.
\begin{Theorem}
Let $p$ be an odd prime and let $k$ be an integer with $1\leq k\leq p-1$.
$$\text{Let}\;\;\Sigma_{p,k}:=\sum_{l=1}^{p-2}\frac{B_l}{l\,k^l}$$
The sum $\Sigma_{p,k}$ is congruent to zero modulo $p$ if and only if the second residue in the $p$-adic expansion of the root of $X^{p-1}+(p-1)!$ of residue $k$ is zero.
\end{Theorem}

\section{Proof of the main theorem and the Fermat quotient $q_p(2)$ modulo $p$}
We start this part by introducing some expressions for various sums involving Bernoulli numbers. These will be useful throughout of the paper.\\
We establish the following lemma which arises from writing Lemma $1$ of the introduction modulo $p^2$.
\begin{Lemma}
\begin{equation}
-p=p!=1+p\Bigg(\frac{B_{p-1}}{1}+\frac{B_{p-3}}{3}+\frac{B_{p-5}}{5}+\dots+\frac{B_2}{p-2}-\frac{3}{2}\Bigg)\mpd
\end{equation}
\begin{equation}
p\,B_{p-1}=-1+p\Bigg(1-2q_p(2)+\sum_{s=1}^{\frac{p-3}{2}}(1-2^{2s+1})\frac{B_{p-(2s+1)}}{2s+1}\Bigg)\;\;\,\,\mpd
\end{equation}
\end{Lemma}
\noindent \textsc{Proof.} The expression for the Stirling numbers modulo $p^2$ is read out of Corollary $2$ of \cite{LEV2}. We have by Eq. $(1)$ of Lemma $1$,
\begin{equation}p!=p\,B_{p-1}-p+\sum_{s=1}^{\frac{p-3}{2}}\frac{p}{2s+1}\,B_{p-2s-1}-\frac{p}{2}+1\mpd\end{equation}
A simple rewriting of $(26)$ is $(24)$.\\
Next, we write Eq. $(2)$ of Lemma $1$ modulo $p^2$. It allows to express the weighted sum of interest in terms of the Agoh-Giuga quotient and of the Fermat quotient. We have, where we used again the expressions for the Stirling numbers modulo $p^2$ taken from Corollary $2$ of \cite{LEV2},
\begin{equation}
-p=(p+1)!=2(p\,B_{p-1}-p)+\sum_{s^{'}=1}^{\frac{p-3}{2}}\frac{p}{2s^{'}+1}\,B_{p-2s^{'}-1}\,2^{2s^{'}+1}-\frac{p}{2}+2^p\mpd
\end{equation}
Congruence $(25)$ is obtained by subtracting $(27)$ and $(26)$. \\

If in the sum of Congruence $(27)$ we make the change of indices corresponding to $2s=p-(2s^{'}+1)$, we get
\begin{equation}
2p\,\sum_{s=1}^{\frac{p-3}{2}}\frac{B_{2s}}{2s\,2^{2s}}=2(p\,B_{p-1}-p)+\frac{p}{2}+2^p\mpd
\end{equation}
We deduce, using also that $B_1=-\frac{1}{2}$ and $B_j=0$ for $j$ odd greater than $1$,
\begin{Lemma}
\begin{equation}
\sum_{k=1}^{p-2}\frac{B_k}{k\,2^k}=\frac{p\,B_{p-1}+1}{p}+q_p(2)-1\mpu
\end{equation}
\end{Lemma}
The next goal is to prove Theorem $3$ which provides a combinatorial interpretation for the Fermat quotient $q_p(2)$ modulo $p$.
To that aim, we will need a result of Lehmer.
\begin{Theorem} (Due to Emma Lehmer \cite{LE}, $1938$)
If $p-1\not|\,2k-2$, then
\begin{equation}
p\,B_{2k}\equiv\frac{1}{2^{2k-1}}\sum_{a=1}^{\frac{p-1}{2}}(p-2a)^{2k}\;\;\text{mod}\,p^3
\end{equation}
\end{Theorem}
Apply Lehmer's result with $2k=p-1$ after noting that $p-1\not|p-3$. \\We get:
\begin{equation}
p\,B_{p-1}=\frac{1}{2^{p-2}}\sum_{a=1}^{\frac{p-1}{2}}(p-2a)^{p-1}\;\;\text{mod}\,p^3
\end{equation}
Expanding the product modulo $p^2$ yields:
\begin{eqnarray}
(p-2a)^{p-1}&=&2^{p-1}a^{p-1}-p(p-1)2^{p-2}a^{p-2}\mpd\\
&=&2^{p-2}(2a^{p-1}+p\,a^{p-2})\qquad\;\;\;\;\,\;\mpd
\end{eqnarray}
It follows that
\begin{equation}
p\,B_{p-1}=2\sum_{a=1}^{\frac{p-1}{2}}a^{p-1}+p\,\sum_{a=1}^{\frac{p-1}{2}}a^{p-2}\mpd
\end{equation}
We stop here to introduce some notations which will be useful throughout the paper.
\begin{Notation}
Denote by $S_k$ the sum of $k$th powers
$$S_k=\sum_{a=1}^{\frac{p-1}{2}}a^{k}$$
\end{Notation}
Using Notation $3$, we summarize congruence $(34)$ in the following lemma.
\begin{Lemma}
\begin{equation}
p\,B_{p-1}=2\,S_{p-1}+p\,H_{\frac{p-1}{2}}\mpd
\end{equation}
\end{Lemma}
\newtheorem{Corollary}{Corollary}
\begin{Remark}
In \cite{SU}, Sun provides a congruence for $H_{\frac{p-1}{2}}$ modulo $p^3$ in terms of the Bernoulli number $B_{p-3}$ and of the Fermat quotient $q_p(2)$. It is the following:
\begin{equation}
H_{\frac{p-1}{2}}=-2q_p(2)+p\,q_p^2(2)-\frac{2}{3}p^2\,q_p^3(2)-\frac{7}{12}\,p^2\,B_{p-3}\;\;\text{mod}\,p^3
\end{equation}
Previously in \cite{LE}, Lehmer had obtained the same congruence modulo $p^2$ only.
\end{Remark}
Another useful result is the following.
\begin{Theorem}
$$S_k=\sum_{r=1}^{\frac{p-1}{2}}r^k=\begin{cases}
(2^{-k+1}-1)B_k\frac{p}{2}\mpd&\text{if $k$ is even}\;\;(I)_k\\
&\\
\big(\frac{1}{2^{k+1}}-1\big)\,\frac{2\,B_{k+1}}{k+1}\mpd&\text{if $k$ is odd}\;\;(II)_k
\end{cases}$$
under the condition $p-1\not|\,k-1$.
\end{Theorem}
The first congruence due to Lehmer for the even $k$'s appears in \cite{LE} while the second congruence for odd $k$'s is due to Mirimanoff and appears in \cite{MIR} at the bottom of page $299$. Mirimanoff's notations for the Bernoulli numbers are the same as Glaisher's notations in \cite{GL}. We provide here a global proof which gathers at once both cases: $k$ is even and $k$ is odd. We will use the Bernoulli polynomials:
\begin{equation}B_n(x)=\sum_{k=0}^n\binom{n}{k}B_{n-k}\,x^k\end{equation}
First and foremost, we have from the Bernoulli formula for the sums of powers:
\begin{equation}
\sum_{r=1}^{\frac{p-1}{2}}r^k=\frac{1}{k+1}\Bigg\lbrace\sum_{l=0}^{k+1}\binom{k+1}{l}B_l\bigg(\frac{p-1}{2}\bigg)^{k+1-l}-B_{k+1}\Bigg\rbrace
\end{equation}
Eq. $(38)$ holds by considering Bernoulli numbers of the second kind, that is $B_1=\frac{1}{2}$ instead of $B_1=-\frac{1}{2}$.
After the change of indices corresponding to $m=k+1-l$, the difference in $(38)$ can be written with Bernoulli polynomials instead, as follows.
\begin{equation}
\sum_{r=1}^{\frac{p-1}{2}}r^k=\frac{B_{k+1}\big(\frac{p-1}{2}\big)-B_{k+1}}{k+1}
\end{equation}
We will use the translation formula
\begin{equation}
B_n(x+y)=\sum_{k=0}^n\binom{n}{k}B_k(x)y^{n-k}
\end{equation}
We apply this formula with $x=-\frac{1}{2}$ and $y=\frac{p}{2}$ and obtain
\begin{equation}
B_{k+1}\big(\frac{p-1}{2}\big)=\sum_{m=0}^{k+1}\binom{k+1}{m}B_m\big(-\frac{1}{2}\big)\frac{p^{k+1-m}}{2^{k+1-m}}
\end{equation}
When working with Bernoulli numbers of the second kind instead of the first kind in the expression for $B_l(-x)$ and conversely in the expression for $B_l(x)$, we have
\begin{equation}
B_l(-x)=(-1)^l\,B_l(x)
\end{equation}
Moreover, we have (see for instance \cite{IR}),
\begin{equation}
B_n\Big(\frac{1}{2}\Big)=\Big(\frac{1}{2^{n-1}}-1\Big)B_n
\end{equation}
Modulo $p^2$, when $p-1\not|\,k-1$ (so that $p$ does not divide the denominator of $B_{k-1}$ by Von Staudt-Clausen's theorem), there are only two terms in the sum of $(41)$, respectively obtained from $m=k+1$ and $m=k$. Using $(42)$ and $(43)$, we get
\begin{equation}
B_{k+1}\big(\frac{p-1}{2}\big)=(-1)^{k+1}\bigg(\frac{1}{2^k}-1\bigg)B_{k+1}+(k+1)(-1)^k\frac{p}{2}\bigg(\frac{1}{2^{k-1}}-1\bigg)B_k
\end{equation}
When $k$ is even (resp odd), the first (resp second) term vanishes. We obtain the congruences of Theorem $17$. \\

We now apply $(I)_{p-1}$. It yields,
\begin{equation}
S_{p-1}=\sum_{r=1}^{\frac{p-1}{2}}r^{p-1}=\bigg(\frac{1}{2^{p-2}}-1\bigg)B_{p-1}\frac{p}{2}\mpd
\end{equation}
We deduce,
\begin{equation}
2^{p-2}S_{p-1}=\bigg(\frac{1}{2}+\frac{1}{2}-2^{p-2}\bigg)B_{p-1}\frac{p}{2}\mpd
\end{equation}
Hence,
\begin{equation}
2^{p-1}S_{p-1}=\frac{p\,B_{p-1}}{2}+\frac{p\,B_{p-1}}{2}(1-2^{p-1})\mpd
\end{equation}
We now derive
\begin{equation}
(1+pq_p(2))S_{p-1}=S_{p-1}-p\,q_p(2)-\frac{p}{2}\frac{(p-1)}{2^{p-1}}\sum_{m=0}^{p-3}(-1)^m\,E(p-2,m)\mpd,
\end{equation}
where we used Congruence $(35)$ of Lemma $4$, Congruence $(36)$ of Remark $2$, and a formula relating the alternating sum of Eulerian numbers to the Bernoulli numbers by
\begin{equation}
(2^{n+1}-1)\frac{B_{n+1}}{n+1}=\frac{1}{2^{n+1}}\sum_{m=0}^{n-1}(-1)^m\,E(n,m),
\end{equation}
applied here with $n+1=p-1$. Modulo $p$, the sum $S_{p-1}$ is simply $-\frac{1}{2}$, hence we further derive
\begin{equation}
\frac{1}{2}p\,q_p(2)=\frac{p}{2}\sum_{m=0}^{p-3}(-1)^m\,E(p-2,m)\mpd
\end{equation}
And after simplifying by $p$ and multiplying by $2$, we obtain
\begin{equation}
q_p(2)=\sum_{m=0}^{p-3}(-1)^m\,E(p-2,m)\mpu
\end{equation}
In order to conclude to Theorem $1$, it suffices to notice that the non alternating sum is the total number of permutations on $p-2$ letters. Moreover, adding the alternating sum with the non alternating sum leads to twice the number of permutations on $p-2$ letters with an even number of ascents. Then $(p-2)!+q_p(2)$ equals this latter number. But $(p-2)!=1\;\text{mod}\,p$. We thus obtain the result of Theorem $1$.

The main Theorem follows from the conjunction of Theorem $1$ and Lemma $3$. As for Theorem $5$, it is derived from the Main Theorem, from the congruence
$$\frac{pB_{p-1}+1}{p}-\frac{1}{p}+\frac{B_{p-1}}{(p-1)2^{p-1}}=1-q_p(2)\;\mpu$$
and from Theorem $1$.

Before closing this part, we say a word about Theorem $6$. Congruence $(9)$ can be derived from Congruence $(24)$; in the main theorem, modulo $p$, the Agoh-Giuga quotient equals the Wilson quotient plus $1$ and still modulo $p$, the second term in the sum is the Fermat quotient minus $1$. We thus obtain $(10)$. Congruence $(11)$ was already discussed before. As for Congruence $(12)$, it comes from a combination of $(24)$, $(25)$, Theorem $3$ and Wolstenholme's theorem. Indeed, Theorem $3$ and Wolstenholme's theorem imply that
$$q_p(2)=\sum_{\begin{array}{l}k=1\\\text{k even}\end{array}}^{p-1}\frac{(-1)}{k}\mpu$$
Then, by $(25)$, we have
$$2\sum_{\begin{array}{l}1\leq k\leq p-2\\\text{k is even}\end{array}}\frac{B_k}{k2^k}=\frac{p\,B_{p-1}+1}{p}+1+\sum_{\begin{array}{l}1\leq k\leq p-2\\\text{k is even}\end{array}}\Big(\frac{B_k-2}{k}\Big)\mpu$$
Finally, by $(24)$, we have
$$\frac{p\,B_{p-1}+1}{p}=\frac{1}{2}+\sum_{\begin{array}{l}1\leq k\leq p-2\\\text{k is even}\end{array}}\frac{B_k}{k}\mpu$$
Congruence $(12)$ follows.

\section{Congruences concerning sums of sums of powers and sums of powers of odd integers}

We begin this part by showing Proposition $1$ which until now remained an independent result of the introduction. Proposition $1$ relates the sum of powers $S_{p-1}$ to the weighted sum studied in the paper. We use a combination of $(25)$ and $(35)$. The terms in $-2q_p(2)$ vanish and we simply get
\begin{equation}
S_{p-1}=\frac{p-1}{2}+\frac{p}{2}\sum_{s=1}^{\frac{p-3}{2}}(1-2^{2s+1})\frac{B_{p-(2s+1)}}{2s+1}\;\text{mod}\,p^2
\end{equation}
It follows that
\begin{equation}
S_{p-1}=\frac{p-1}{2}+\frac{1}{2}(\frac{p}{2}-1-p\,B_{p-1})-p\sum_{s=1}^{\frac{p-3}{2}}2^{2s}\frac{B_{p-2s-1}}{2s+1}\mpd,
\end{equation}
where we used $(24)$. \\
From there we easily derive the congruence of Proposition $1$. \\Note with $p\,B_{p-1}=-1\;\text{mod}\,p$ that the formula is correct modulo $p$. \\
We can now derive the new Proposition $2$ below.
\begin{Proposition}
\begin{equation}\sum_{\begin{array}{l}1\leq a\leq p-1\\\text{a odd}\end{array}}a^{p-1}=\frac{1}{2}+p\,B_{p-1}-\frac{3p}{8}-\frac{p}{4}\sum_{k=1}^{\frac{p-3}{2}}\frac{B_{2k}}{k2^k}\mpd
\end{equation}
\end{Proposition}
\textsc{Proof.} We have,
\begin{eqnarray}
\sum_{\begin{array}{l}1\leq a\leq p-1\\\;\,\text{a odd}\end{array}}a^{p-1}&=&\sum_{a=1}^{p-1}a^{p-1}-\sum_{a=1}^{\frac{p-1}{2}}(2a)^{p-1}\\
&=&p\,B_{p-1}+2^{p-1}+2^{p-2}p\,B_{p-1}-\frac{p}{2}\bigg(\frac{3}{2}+\sum_{k=1}^{\frac{p-3}{2}}\frac{B_{2k}}{k2^{2k}}\bigg)\notag\\&&\qquad\qquad\qquad\qquad
\qquad\qquad\qquad\qquad\mpd
\end{eqnarray}
For the sum of $(p-1)$-th powers of the first $p-1$ integers modulo $p^2$, the reader is referred to \cite{LEV2} or to Congruence $(5.1)$ of \cite{SU}.
Next, we have
\begin{equation}
2^{p-2}=\frac{1}{2}+\frac{1}{2}p\,q_p(2),
\end{equation}
and so
\begin{equation}
2^{p-2}\,p\,B_{p-1}=\frac{1}{2}p\,B_{p-1}-\frac{1}{2}p\,q_p(2)\mpd
\end{equation}
Before moving forward in the computation, we introduce a new notation.
\begin{Notation}
$$\mathcal{S}_p:=\sum_{k=1}^{\frac{p-3}{2}}\frac{B_{2k}}{k2^{2k}}$$
\end{Notation}
Plugging $(57)$ and $(58)$ into $(56)$ yields after simplifying,
\begin{equation}
\sum_{\begin{array}{l}1\leq a\leq p-1\\\;\,\text{a odd}\end{array}}a^{p-1}=\frac{3}{2}\,pB_{p-1}+1+\frac{1}{2}pq_p(2)-\frac{p}{2}(\frac{3}{2}+\mathcal{S}_p)\mpd
\end{equation}
Now, from Congruence $(25)$ we derive:
\begin{equation}
p\,B_{p-1}-p+1=-2p\,q_p(2)+\frac{p}{2}-p\,B_{p-1}-1+p\,\mathcal{S}_p\mpd
\end{equation}
From there, we get
\begin{equation}
\frac{1}{2}pq_p(2)=-\frac{1}{2}-\frac{1}{2}p\,B_{p-1}+\frac{3p}{8}+\frac{1}{4}p\,\mathcal{S}_p\mpd
\end{equation}
And so, plugging back into $(59)$, we obtain Congruence $(54)$ of Proposition $2$. \\

More generally, let $k$ be an integer with $1<k\leq \frac{p-1}{2}$. We can express the sum
$$1^{2k-1}+3^{2k-1}+5^{2k-1}+\dots +(p-2)^{2k-1}$$
in terms of an alternating sum of Eulerian numbers modulo $p^2$. Indeed, we have:
\begin{eqnarray}
\sum_{\begin{array}{l}x=1\\\text{x odd}\end{array}}^{p-1}x^{2k-1}&=&\sum_{x=1}^{p-1}x^{2k-1}-2^{2k-1}\sum_{x=1}^{\frac{p-1}{2}}x^{2k-1}\\
&=&-2^{2k-1}\sum_{x=1}^{\frac{p-1}{2}}x^{2k-1}\;\;\;\;\;\;\;\;\mpd\\
&=&\big(2^{2k-1}-\frac{2^{2k-1}}{2^{2k}}\big)\frac{2B_{2k}}{2k}\mpd\\
&=&(2^{2k}-1)\frac{B_{2k}}{2k}\;\;\;\;\qquad\;\;\;\mpd
\end{eqnarray}
The first sum in the right hand side of $(62)$ is congruent to zero modulo $p^2$. This fact gets extensively discussed in \cite{LEV2}.
Congruence $(64)$ is then obtained by applying $(II)_{2k-1}$ of Theorem $17$ since $k$ was chosen so that $p-1\not| 2k-2$. An application of $(49)$ now yields the result stated in the following proposition.
\begin{Proposition}
Let $k$ be an integer with $1<k\leq \frac{p-1}{2}$. Then, we have
\begin{equation}
\sum_{\begin{array}{l}x=1\\\text{x odd}\end{array}}^{p-1}x^{2k-1}=\frac{1}{2^{2k}}\sum_{m=0}^{2k-2}(-1)^m\,E(2k-1,m)\mpd
\end{equation}
\end{Proposition}
\noindent In particular, applying $(66)$ with $k=\frac{p-1}{2}$ yields:
\begin{equation}
\sum_{\begin{array}{l}x=1\\\text{x odd}\end{array}}^{p-1}x^{p-2}=\frac{1}{2^{p-1}}\sum_{m=0}^{p-3}(-1)^m\,E(p-2,m)\mpd
\end{equation}
It follows that, where we used Congruence $(51)$,
\begin{equation}
(1+pq_p(2))\,\sum_{\begin{array}{l}x=1\\\text{x odd}\end{array}}^{p-1}x^{p-2}=(q_p(2))_0+p\,\Bigg(\sum_{m=0}^{p-3}(-1)^m\,E(p-2,m)\Bigg)_1\mpd
\end{equation}
Moreover, modulo $p$, the sum we are computing is nothing else than $H^{'}_{p-1}$ which by Theorem $3$ is also $q_p(2)$ modulo $p$.
We thus derive:
\begin{equation}
\begin{split}\big(q_p(2)\big)_0+p\,\Bigg\lbrace \Big(\sum_{\begin{array}{l}x=1\\\text{x odd}\end{array}}^{p-1}x^{p-2}\Big)_1&+\big(q_p(2)\big)_0^2\Bigg\rbrace=\big(q_p(2)\big)_0\\&+p\,\Bigg(\sum_{m=0}^{p-3}(-1)^m\,E(p-2,m)\Bigg)_1\mpd
\end{split}\end{equation}
It comes:
\begin{equation}
\Big(\sum_{\begin{array}{l}x=1\\\text{x odd}\end{array}}^{p-1}x^{p-2}\Big)_1=\Bigg(\sum_{m=0}^{p-3}(-1)^m\,E(p-2,m)\Bigg)_1-\big(q_p(2)\big)_0^2\;\;\;\;\,\mpu
\end{equation}
Moreover, we have
\begin{equation}
(p-2)!+\sum_{m=0}^{p-3}(-1)^m\,E(p-2,m)=2\,\sum_{k=0}^{\frac{p-3}{2}}E(p-2,2k)
\end{equation}
And so,
\begin{equation}\begin{split}
-(p+1)(p\,B_{p-1}-p)&+(q_p(2))_0+\\&p\Bigg(\sum_{m=0}^{p-3}(-1)^m\,E(p-2,m)\Bigg)_1=2\,N_{p-2}\mpd
\end{split}\end{equation}
Thus,
\begin{equation}\begin{split}
2p-p\,B_{p-1}&+(q_p(2))_0+\\&p\Bigg(\sum_{m=0}^{p-3}(-1)^m\,E(p-2,m)\Bigg)_1=(2N_{p-2})_0+p(2N_{p-2})_1\mpd
\end{split}\end{equation}
From there, we retrieve in particular the fact that $$(q_p(2))_0=(2N_{p-2})_0-1\,\mpu,$$ which is the result of Theorem $1$.
In order to proceed further, we will need a simple lemma whose statement appears right below.
\newtheorem*{Lem}{Lemma}
\begin{Lem} We have using Notation $2$:
\begin{equation}\Big[(q_p(2))_0+1-(2N_{p-2})_0\Big]_1=-\big((2N_{p-2})_0-1\big)_1\mpu\end{equation}
\end{Lem}
\textsc{Proof of the lemma}. Indeed, we have:
\begin{eqnarray*}
\Big[(q_p(2))_0+1-(2N_{p-2})_0\Big]_1&=&\Bigg[-\Big((2N_{p-2})_0-1-\big((2N_{p-2})_0-1\big)_0\Big)\Bigg]_1\\
&=&-\Big((2N_{p-2})_0-1-\big((2N_{p-2})_0-1\big)_0\Big)_1\negthickspace\!\qquad\mpu\\
&=&-\big((2N_{p-2})_0-1\big)_1\,\qquad\qquad\qquad\qquad\qquad\mpu
\end{eqnarray*}
By combining Congruences $(70)$ and $(73)$ and using the lemma, we obtain:
\begin{equation}
\Big(\sum_{\begin{array}{l}x=1\\\text{x odd}\end{array}}^{p-1}x^{p-2}\Big)_1=(p\,B_{p-1})_1+(2N_{p-2})_1+\big((2N_{p-2})_0-1\big)_1-\big((2N_{p-2})_0-1\big)^2-2
\,\mpu
\end{equation}
From there, point $(i)$ of Theorem $8$ follows directly. For point $(ii)$, we cannot use Proposition $3$, hence we make a direct calculation.
We have,
\begin{eqnarray}
\sum_{\begin{array}{l}x=1\\\text{x odd}\end{array}}^{p-1}x^{p-1}&=&\sum_{x=1}^{p-1}x^{p-1}-2^{p-1}\sum_{x=1}^{\frac{p-1}{2}}x^{p-1}\\
&=&\frac{1}{2}(p\,B_{p-1}-pq_p(2))\qquad\qquad\qquad\qquad\;\;\;\;\,\,\,\mpd\\
&=&-\frac{1}{2}+\frac{p}{2}\bigg((p\,B_{p-1})_1-(2N_{p-2})_0+1\bigg)\qquad\mpd
\end{eqnarray}
Congruence $(77)$ is obtained by applying $(I)_{p-1}$ and $(78)$ follows from applying Theorem $1$. Theorem $8$ is thus entirely proven. \\

We now deal with sums of sums of powers. \\

\noindent Using Notation $3$, the sum of the first $\frac{p-1}{2}$ odd powers of the first $\frac{p-1}{2}$ integers is
$$S_1+S_3+\dots+S_{p-4}+S_{p-2}$$
First we sum $T:=S_1+S_3+\dots+S_{p-4}$, using $(II)_{2k-1}$ of Theorem $17$. It comes:
\begin{equation}
p\sum_{k=1}^{\frac{p-3}{2}}S_{2k-1}=p\,\sum_{k=1}^{\frac{p-3}{2}}\frac{B_{2k}}{k2^{2k}}-2p\sum_{k=1}^{\frac{p-3}{2}}\frac{B_{2k}}{2k}\mpd
\end{equation}
Since we work modulo $p^2$, the second sum may be replaced with $$-\sum_{k=1}^{\frac{p-3}{2}}\frac{B_{2k}}{p-2k}$$ instead, and so we get
by using $(24)$ together with Proposition $1$,
\begin{eqnarray}
p\sum_{k=1}^{\frac{p-3}{2}}S_{2k-1}&=&2\bigg(S_{p-1}+1+\frac{1}{2}p\,B_{p-1}\bigg)-\frac{3}{2}p+2\bigg(\frac{3}{2}p-p\,B_{p-1}-1-p\bigg)\mpd\notag\\
&=&-\frac{1}{2}p-p\,B_{p-1}+2S_{p-1}\mpd\\
&=&-\frac{1}{2}p+2p\,q_p(2)\;\;\;\mpd
\end{eqnarray}
The last congruence is obtained by using $(35)$ and $(36)$ which jointly provide
\begin{equation}2S_{p-1}=p\,B_{p-1}+2pq_p(2)\;\;\;\text{mod}\,p^2\end{equation}
Then, after simplifying by $p$, we get:
\begin{equation}
T=-\frac{1}{2}+2q_p(2)\mpu
\end{equation}
Further,
\begin{equation}
S_{p-2}=H_{\frac{p-1}{2}}=-2q_p(2)\mpu
\end{equation}
So, adding all the terms in the considered sum now yields
\begin{equation}
S_1+S_3+\dots +S_{p-4}+S_{p-2}=-\frac{1}{2}=S_{p-1}\mpu,
\end{equation}
as announced in the statement $(i)$ of Theorem $7$.\\ Point $(ii)$ then simply follows from the fact that the sum of any odd power less than or equal to $p-2$ of the first $p-1$ integers is divisible by $p$. Hence, so is the sum over the first $\frac{p-1}{2}$ odd powers.

\section{The Fermat quotient $q_p(2)$ modulo $p^2$.}
In this part, we show how our work and Emma Lehmer's work combined allow to get to the next $p$-power in the $p$-adic expansion of the Fermat quotient $q_p(2)$. The congruence which we will use and which is Congruence $(34)$ of \cite{LE} is the following. Lehmer's notation for the Fermat quotient $q_p(2)$ is simply $q_2$ and her $w_p$ denotes the Wilson quotient.
\begin{Theorem} (Emma Lehmer \cite{LE}, $1938$)
\begin{equation}
\sum_{r=1}^{\frac{p-1}{2}}r^{p-2}=-2q_2(1-pw_p)+2pq_2^2\mpd
\end{equation}
\end{Theorem}
Since $$\sum_{x=1}^{p-1}x^{p-2}=p\,B_{p-2}=0\mpd,$$
we have:
\begin{eqnarray}
\sum_{\begin{array}{l}x=1\\\text{x odd}\end{array}}^{p-1}x^{p-2}&=&-\sum_{x=1}^{\frac{p-1}{2}}2^{p-2}x^{p-2}\qquad\qquad\qquad\qquad\qquad\;\;\;\;\,\,\mpd\\
&=&-\frac{1}{2}(1+p\,q_p(2))(-2q_p(2)(1-p\,w_p)+2p\,q_p(2)^2)\notag\\&&\qquad\qquad\qquad\qquad\qquad\qquad\qquad\qquad\qquad\mpd\\
&=&q_p(2)(1-p\,w_p)\qquad\qquad\qquad\qquad\qquad\;\;\;\;\,\,\mpd
\end{eqnarray}
By comparing Congruence $(89)$ with Congruence $(15)$ of Theorem $8$ (point $(i)$), 
we conclude to Congruence $(17)$ of Theorem $9$. And since by Remark $1$, we also have (cf Glaisher's result),
\begin{equation}
q_p(2)=H^{'}_{p-1}+\frac{1}{2}p(H^{'}_{p-1})^2\mpd,
\end{equation}
then Congruences $(18)$ and $(19)$ of Theorem $10$ follow.

We end this part by showing that $(N_{p-2})_0$ can be expressed in terms of a sum of harmonic numbers with odd indices (cf Proposition $4$ below) and we shall also get an expression to the next power of $p$ (cf forthcoming Proposition $8$)
\begin{Proposition}
$$N_{p-2}=\sum_{\begin{array}{l}m=1\\\text{m odd}\end{array}}^{p-2}H_m\mpu$$
\end{Proposition}
\textsc{Proof.} Let $m$ be an integer with $2m\leq p-3$. A classical formula on Eulerian numbers reads
\begin{equation}
E(p-2,2m)=\sum_{k=0}^{2m}(-1)^k\binom{p-1}{k}(2m+1-k)^{p-2}
\end{equation}
\begin{Lemma} The following congruence holds:
$$\forall 1\leq k\leq p-1,\,\binom{p-1}{k}= (-1)^k\mpu$$
\end{Lemma}
\textsc{Proof of the lemma.} As is well known, $$\binom{p}{k}=\binom{p-1}{k-1}+\binom{p-1}{k}$$ and
$\binom{p}{k}$ is divisible by $p$ for each integer $k$ with $1\leq k\leq p-1$. It follows that
$$\binom{p-1}{k}=-\binom{p-1}{k-1}\mpu$$
Then, working inductively, we obtain the result of the lemma. We note that conversely, it was recently proven by Mestrovic in \cite{ME} that if $n>1$ and $q>1$ are integers such that
$$\forall 0\leq k\leq n-1,\,\binom{n-1}{k}=(-1)^k\,\text{mod}\,q,$$
then $q$ is a prime and $n$ is a power of $q$. \\

\noindent Applying Lemma $5$ to Eq. $(91)$ yields:
\begin{eqnarray}
E(p-2,2m)&=&\sum_{k=0}^{2m}\;(2m+1-k)^{p-2}\mpu\notag\\&&\\
&=&\sum_{k=0}^{2m}\;\frac{1}{2m+1-k}\qquad\;\mpu\\
&&\notag\\
&=&H_{2m+1}\qquad\qquad\qquad\mpu
\end{eqnarray}
The result of Proposition $4$ follows. \\

\noindent By the same arguments as before, we have for each $m$ with $0\leq m\leq\frac{p-5}{2}$,
$$E(p-2,2m+1)=H_{2m+2}\mpu$$
It follows that
\begin{eqnarray}
(p-2)!=\sum_{k=0}^{p-3}E(p-2,k)&=&\sum_{m=0}^{\frac{p-3}{2}}E(p-2,2m)+\sum_{m=0}^{\frac{p-5}{2}}E(p-2,2m+1)\notag\\
&&\\
&=&\sum_{\begin{array}{l}m=1\\\text{m odd}\end{array}}^{p-2}H_m+H_2+H_4+\dots+H_{p-3}\notag\\&&\qquad\qquad\qquad\qquad\qquad\qquad\qquad\mpu
\end{eqnarray}
And since by Wolstenholme's theorem we also have $H_{p-1}=0\,\text{mod}\,p$, we conclude to the congruence:
\begin{Proposition}
\begin{equation}
\sum_{k=1}^{p-1}H_k=1\mpu
\end{equation}
\end{Proposition}

\noindent This is Corollary $1.5$ of \cite{MES}, which we retrieved here by a different method. \\
In \cite{MES}, Mestrovic and Andjic prove as a special case of a more general result the stronger statement.
\begin{Theorem}(Mestrovic and Andjic \cite{MES}, $2017$)
$$\sum_{k=1}^{p-1}H_k=1-p\;\text{mod}\,p^3$$
\end{Theorem}
\noindent Proposition $4$ joint with Eq. $(18)$ of Theorem $10$ imply the unexpected result of Theorem $11$ from the introduction. \\
We now investigate $N_{p-2}$ modulo $p^2$. To that aim, we first show the following proposition.
\begin{Proposition}
If $p\geq 3$ is a prime, then
\begin{equation}
\binom{p-1}{k}=(-1)^k(1-p\,H_k)\mpd
\end{equation}
\end{Proposition}
Indeed, the result follows from the following two lemmas.
\begin{Lemma}
$$\binom{p-1}{k}=\binom{p}{k}-\binom{p}{k-1}+\binom{p}{k-2}-\dots+(-1)^{k-1}\binom{p}{1}+(-1)^k$$
\end{Lemma}
\begin{Lemma}
$$\binom{p}{k}=p\;\frac{(-1)^{k+1}}{k}\mpd$$
\end{Lemma}
\noindent We note that Proposition $6$ is part of a stronger congruence by Zhi-Wei Sun.
\begin{Theorem}(Z-W. Sun, Lemma $2.1$ $(2.2)$ of \cite{SU3}, $2012$)
If $p\geq 3$ is a prime, then
\begin{equation}
\binom{p-1}{k}=(-1)^k(1-p\,H_k+\frac{p^2}{2}(H_k^2-H_{k,2}))\;\text{mod}\,p^3
\end{equation}
\end{Theorem}
We deduce an expression for $E(p-2,2m)$ modulo $p^2$.
\begin{Proposition}
\begin{equation}
E(p-2,2m)=\sum_{s=1}^{2m+1}s^{p-2}-p\;\sum_{K=p-(2m+1)}^{p-2}\frac{H_K}{K+(2m+2)}\mpd
\end{equation}
\end{Proposition}
\textsc{Proof.} We have by applying Proposition $6$ in Eq. $(91)$,
\begin{equation}
E(p-2,2m)=\sum_{s=1}^{2m+1}s^{p-2}-p\,\sum_{s=1}^{2m}\frac{H_{2m+1-s}}{s}\mpd
\end{equation}
Moreover, we have by using Wolstenholme's theorem,
\begin{equation}
H_{p-(2m+1-s+1)}=H_{2m+1-s}\mpu
\end{equation}
Indeed, we proceed by induction. We have,
\begin{eqnarray*}H_{p-2}&=&H_{p-1}-\frac{1}{p-1}=1=H_1\mpu\\
H_{p-3}&=&H_{p-2}-\frac{1}{p-2}=1+\frac{1}{2}=H_2\mpu\\
&\vdots&\\
H_{p-k}&=&H_{p-k+1}-\frac{1}{p-k+1}=1+\frac{1}{2}+\dots+\frac{1}{k-1}=H_{k-1}\mpu
\end{eqnarray*}
Congruence $(102)$ follows immediately. We derive Proposition $8$ below.
\begin{Proposition}
\begin{equation}
N_{p-2}=\sum_{m=0}^{\frac{p-3}{2}}\Bigg(S_{2m+1,\,p-2}-p\;\sum_{K=p-(2m+1)}^{p-2}\frac{H_K}{K+(2m+2)}\Bigg)\mpd
\end{equation}
\end{Proposition}

\noindent We note that our formula is consistent with our Proposition $4$ modulo $p$.

Proposition $8$ is only the starting point of further investigations needed. Namely, expressing Proposition $4$ and $8$ in terms of Bernoulli numbers requests more work. It is to expect that working out further the sums in Proposition $8$ could lead to conclusions on the existence or non-existence of Wilson and super Wieferich primes.

\section{Ending notifications and proof of Theorem $15$}

When we were about to submit this paper, we learned that Zhi-Wei Sun, Zhi-Hong Sun's brother, and his co-author Li-Lu Zhao studied congruences concerning the same weighted sums, where the (non divided) Bernoulli numbers are replaced with harmonic numbers or generalized harmonic numbers instead. Their result taken from \cite{SUB} is the following, where notations are standard.

\begin{Theorem} (Zhi-Wei Sun and Li-Lu Zhao \cite{SUB}, $2013$)
Let $p$ be a prime with $p>3$. Then,
\begin{eqnarray}
\sum_{k=1}^{p-1}\frac{H_k}{k\,2^k}&=&\frac{7}{24}\,p\,B_{p-3}\;\;(mod\,p^2)\\
\sum_{k=1}^{p-1}\frac{H_{k,2}}{k2^k}&=&-\frac{3}{8}\,B_{p-3}\;\;(mod\,p)
\end{eqnarray}
\end{Theorem}
Previously, Zhi-Wei Sun in \cite{SU3} had shown under the same conditions that $$\sum_{k=1}^{p-1}\frac{H_k}{k2^k}=0\;\;\text{mod}\,p$$
We note that if $(p,p-3)$ is an irregular pair, the powers of two weighted sum of the first $(p-1)$ divided harmonic numbers (resp generalized harmonic numbers) is divisible by $p^2$ (resp by $p$). Such a prime is called a Wolstenholme prime as for these primes, Wolstenholme's theorem \cite{WO1} holds to the next power of $p$, that is $H_{p-1}=0\;\text{mod}\,p^3$ and $H_{p-1,2}=0\;\text{mod}\,p^2$.\\

Here, our congruence $(10)$ implies that
\begin{equation}\sum_{k=1}^{p-2}\frac{B_k}{k2^k}=0\;\text{mod}\,p\Leftrightarrow q_p(2)=-w_p\;\text{mod}\,p\end{equation}

We provide a polynomial interpretation of Congruence $(106)$. Our polynomial actor is $X^{p-1}+(p-1)!$, whose factored form in the ring of $p$-adic integers appears in the proposition below.

\begin{Proposition} We have the factorization in $\mathbb{Z}_p[X]$ and subsequent equality
\begin{eqnarray}X^{p-1}+(p-1)!&=&(X-1-pt_1)(X-2-pt_2)\dots (X-(p-1)-pt_{p-1})\notag\\
q_p(2)+w_p&=&-t_2(1-pt_1)\Big[\prod_{i=1}^{p-3}(i+p\,t_{i+2})\Big]
\end{eqnarray}
\end{Proposition}

\textsc{Proof.} By \cite{LEV2}, we may factor the given polynomial in $\mathbb{Z}_p[X]$ as in the statement, where the $t_i$'s with $1\leq i\leq p-1$ are some adequate $p$-adic integers. Specializing $X=2$ leads to identity $(107)$.\\

\noindent We retrieve from $(107)$ the fact that $t_2^{(0)}=2(q_p(2)+w_p)\,\mpu$. This is a special case of Lemma $1$ of \cite{LEV2}.
Thus, we see that Congruence $(106)$ holds if and only if the second residue in the $p$-adic expansion of the root of $X^{p-1}+(p-1)!$ of residue $2$ is zero.\\
\indent In fact, we can generalize this result to powers of $k$ weighted sums of divided Bernoulli numbers, when $k$ is an integer with $3\leq k\leq p-1$. It suffices to establish the following two lemmas.
\begin{Lemma}
Let $k$ be an integer with $1\leq k\leq p-1$. Then,
$$p=-\sum_{s=1}^{p}\left[\begin{array}{l}p\\s\end{array}\right]\,k^s\mpd$$
\end{Lemma}
\begin{Lemma}
Let $k$ be an integer with $1\leq k\leq p-1$. Let $q_p(k)$ denote the Fermat quotient to base $k$, with standard notations. Then,
$$\sum_{l=1}^{p-2}\frac{B_l}{l\,k^l}=0\mpu\Leftrightarrow q_p(k)=-w_p\mpu$$
\end{Lemma}
\textsc{Proof.} Lemma $8$ follows from specializing $x=-k$ in Eq. $(4)$ and applying Wilson's theorem modulo $p$. First, we have from the specialization:
$$\frac{(p-1+k)!}{(k-1)!}=\sum_{s=1}^p\left[\begin{array}{l}p\\s\end{array}\right]k^s$$
Then, we write
\begin{eqnarray}
(p-1+k)!&=&(p-1)!p(p+1)\dots(p-1+k)\\
&=&-p\,(k-1)!\mpd
\end{eqnarray}

\noindent Then, Lemma $9$ can be derived from Lemma $8$ by using the expressions for the Stirling numbers modulo $p^2$ as provided for instance in \cite{LEV2}. The calculations are similar to those realized before and are left to the reader. \\
Moreover, we have by using the specialization $X=k$ in both factored and expanded forms of the polynomial $X^{p-1}+(p-1)!$,
\begin{equation}
q_p(k)+w_p=(-1)^{k+1}\,t_k^{(0)}\,(k-1)!\,(p-1-k)!\mpu
\end{equation}
This implies by the Gauss lemma that $t_k^{(0)}=0$. Hence the result of Theorem $15$ that was originally announced at the end of the introduction. Note that the theorem holds even when $k=1$.\\

\textsc{Email address:} \textit{clairelevaillant@yahoo.fr}


\begin{thebibliography}{ll}

\bibitem{BE} N.G.W.H. Beeger, On a new case of the congruence $2^{p-1}\equiv 1\;(p^2)$, Messenger of Mathematics $51$ $(1922)$ $149-150$
\bibitem{CDP} R. Crandall, K. Dilcher and C. Pomerance, A search for Wieferich and Wilson primes, Math. Comput. Vol. $66$, No. $217$ $(1997)$ $433-449$
\bibitem{DO} F.G. Dorais, D. Klyve, A Wieferich prime search up to $6.7\times\,10^{15}$, J. Integer Sequences $14$ (9) $(2011)$
\bibitem{EI} G. Eisenstein, Eine neue Gattung zahlentheoreticher Funktionen, welche von zwei Elementen abh\"angen und durch gewisse lineare Funktional-Gleichungen definirt werden, Berichte K\"onigl. Preu\ss$\;\,$ Akad. Wiss. Berlin $15$ $(1850)$ $36-42$
\bibitem{JF} J. Faulhaber, Academia Algebrae, Darinnen die miraculosische Inventiones zu den h\"ochsten Cossen weiters continuirt und prolifiert werden $1631$, QA$154.8\,F3\,1631afMATH$ at Stanford University Libraries by courtesy of Donald E. Knuth.
\bibitem{GL} J.W.L. Glaisher, On the residues of the sums of products of the first $p-1$ numbers and their powers, to modulus $p^2$ or $p^3$, Quarterly J. Math. $31$ $(1900)$ $321-353$
\bibitem{GL2} J.W.L. Glaisher, On the residues of $r^{p-1}$ to modulus $p^2$, $p^3$, etc., Q.J. Math. Oxford $32$ $1900-1901$ $1-27$
\bibitem{IR} K. Ireland and M. Rosen, A classical introduction to modern number theory, Springer NY $1982$
\bibitem{LE} E. Lehmer, On congruences involving Bernoulli numbers and the quotients of Fermat and Wilson, Ann. Math. $39$ $(1938)$ $350-360$
\bibitem{LEV1} C. Levaillant, A quantum combinatorial approach for computing a tetrahedral network of Jones-Wenzl projectors, arXiv:$1301.1733$
\bibitem{LEV2} C. Levaillant, Wilson's theorem modulo $p^2$ derived from Faulhaber polynomials, arXiv:$1912.06652$
\bibitem{MI} W. Meissner, \"Uber die Teilbarkeit von $2^p-2$ durch das Quadrat der Primzhal $p=1093$, Sitzungsber. Akad. d. Wiss. Berlin, $(1913)$ $663-667$
\bibitem{ME} R. Mestrovic, A primality criterion based on a Lucas'congruence, arXiv:$1407.7894v1$
\bibitem{MES} R. Mestrovic and M. Andjic, Certain congruences for harmonic numbers, Mathematica Montisnigri, Vol. XXXVIII $(2017)$
\bibitem{MIR} D. Mirimanoff, Sur la congruence $(r^{p-1}-1):p\equiv q_r\;(mod.p)$, J. Reine und Angew. Math. $115$ $(1895)$ $295-300$
\bibitem{SI} J. Silverman, Wieferich's criterion, J. Number Theory $30$ $(1988)$ $226-237$
\bibitem{ST} M. Stern, Einige Bermerkungen \"uber die Congruenz $\frac{(r^p-r)}{p}\equiv a\;(mod\,p)$, J. Reine und Angew. Math. $100$ $(1887)$ $182-188$
\bibitem{SU} Z-H. Sun, Congruences concerning Bernoulli numbers and Bernoulli polynomials, Discrete Applied Math. $105$ $(2000)$ $193-223$
\bibitem{SU3} Z-W. Sun, Arithmetic theory of harmonic numbers, Proc. Amer. Math. Soc. $140$ $(2012)$ $415-428$
\bibitem{SUB} Z-W. Sun and L.L. Zhao, Arithmetic theory of harmonic numbers $(II)$, Colloq. Math. $130$ $(2013)$, no. $1$ $67-78$
\bibitem{VA} H.S. Vandiver, Symmetric functions formed by systems of elements of a finite algebra and their connection with Fermat's quotient and Bernoulli's numbers, Ann. Math. $18$ $(1917)$ $105-114$
\bibitem{WI} A. Wieferich, Zum letzten Fermat'schen Theorem, J. Reine und Angew. Math. $136$ $(1909)$ $293-302$
\bibitem{WO1} J. Wolstenholme, On certain properties of prime numbers, Quaterly J. of Pure and Applied Math. Vol $5$ $(1862)$ $35-39$
\bibitem{WO2} J. Worpitzky, Studien \"uber die Bernoullischen und Eulerschen Zahlen, J. Reine und Angew. Math. $(1883)$ $94$ $203-232$

\end{thebibliography}
\end{document}